\documentclass[smallextended]{article}
\usepackage{amsfonts}
\usepackage{graphicx}
\usepackage{amsmath}
\usepackage{url}
\newtheorem{theorem}{Theorem}

\newtheorem{corollary}{Corollary}

\newtheorem{definition}{Definition}

\newtheorem{lemma}{Lemma}

\newtheorem{property}{Property}

\newcommand{\be}{\begin{equation}}
\newcommand{\ee}{\end{equation}}
\input{tcilatex}
\begin{document}
\title{\textbf {Ramanujan Summation and the Exponential Generating Function $
\displaystyle{\sum_{k=0}^{\infty }\frac{z^{k}}{k!}\zeta ^{\prime }(-k)}$}}
\author{B. Candelpergher$^{1}$, H. Gopalkrishna Gadiyar$^{2}$ and R. Padma$^{2}$\\
~~\\
$^{1}$Laboratory J.A. Dieudonn\'{e}, UMRCNRS No. 6621\\
University of Nice-Sophia Antipolis\\
Parc Valrose, 06108 Nice Cedex 2, France\\
E-mail: Bernard.CANDELPERGHER@unice.fr \\
~~\\
$^2$AU-KBC Research Centre, M. I. T. Campus of \\Anna University, Chromepet, Chennai 600 044, INDIA\\E-mail: \{gadiyar, padma\}@au-kbc.org}
\date{~~}
\maketitle

\begin{abstract}
In the sixth chapter of his notebooks Ramanujan introduced a method of summing divergent series which assigns to the series the value of the associated Euler-MacLaurin constant that arises by applying the Euler-MacLaurin summation formula to the partial sums of the series. This method is now called the Ramanujan summation process. In this paper we calculate the Ramanujan sum of the exponential generating functions $\sum_{n\geq 1}\log n~e^{nz}$ and $\sum_{n\geq 1}H_n^{(j)}~e^{-nz}$ where $H_n^{(j)}=\sum_{m=1}^n \frac{1}{m^j}$. We find a surprising relation between the two sums when $j=1$ from which follows a formula that connects the derivatives of the Riemann zeta - function at the negative integers to the Ramanujan summation of the divergent Euler sums $\sum_{n\ge 1} n^kH_n,~k~\ge ~0$, where $H_n= H_n^{(1)}$. Further, we express our results on the Ramanujan summation in terms of the classical summation process called the Borel sum.
\end{abstract}


{\bf Mathematics Subject Classfication (2000)} 11M06, 65B15, 40G99

\section{Introduction} Let
\be
\frac{ze^{xz}}{e^z-1} = \sum_{n=0}^\infty B_n(x) \frac{z^n}{n!} \, . \label{eq:bernoulli}
\ee
$B_n(x)$ is the $n^{th}$ Bernoulli polynomial and $B_n(0) = B_n$ is the  $n^{th}$ Bernoulli number. At the beginning of the sixth chapter of his Notebooks \cite{berndt}, Ramanujan writes the Euler - MacLaurin formula for the partial sums 
\begin{equation}
a(1)+a(2)+\cdot +a(x-1) = C+\int_1^x a(t) dt + \sum_{k \ge 1} \frac{B_k}{k!} \partial ^{k-1}(a(x))
\end{equation}
of an infinite series ${\displaystyle \sum_{n = 1}^\infty  a(n)}$  and assigns the value of the constant $`C\textrm{'}$  the sum of the series. For example, one has 
\begin{equation}
\sum_{n \ge 1}^{\cal R} \frac{1}{n} = \gamma
\end{equation}
by this process, where ${\cal R}$ denotes the Ramanujan summation and $\gamma $ is the Euler - Mascheroni constant This follows from the formula \cite{apostol}
\be
\sum_{n \le x} \frac{1}{n} = \log x +\gamma +O\left (\frac{1}{x}\right ) \, .
\ee
It was pointed out by Hardy (Section 6, Chapter 13 of \cite{hardy}) that this definition makes the Ramanujan sum of a series ambiguous as the sum depends on a particular parameter that appears as a lower limit of an integral. In \cite{Candelpergher:Coppo:Delabaere}, Candelpergher, Coppo and Delabaere formulated the Ramanujan summation in a rigorous manner so that the uniqueness of the sum of the series is established. They also proved several properties of the Ramanujan summation and summed many well known divergent and convergent series. One of the results they obtained using this summation process was that
\begin{equation}
\sum_{n \ge 1}^{\cal R}H_n = \frac{3}{2}\gamma +\frac{1}{2}-\frac{1}{2}\log (2\pi) \, , \label{eq:candel}
\end{equation}
where $\displaystyle{H_n=\sum_{j=1}^{n}\frac{1}{j}}$ is the $n^{th}$ harmonic number and $H_0=0$. Note that 
\be
\zeta '(0) = -\frac{1}{2} \log 2 \pi \, .
\ee
So the natural question to ask now is whether $\displaystyle{\sum_{n \ge 1}^{\cal R}n^k H_n }$ is related to $\zeta'(-k)$ for $k \ge 1$. 

A heuristic argument suggesting an identity involving $\displaystyle{\sum_{n = 1}^{\infty}n^k H_n }$ and $\zeta'(-k)$ is given in \cite{gadiyar:padma}. Also, recall the Euler sum identities \cite{Flajolet:Salvy}
\begin{eqnarray}
\sum_{n=1}^\infty \frac{H_n}{n^2}&=&2 \zeta (3)\\
\sum_{n=1}^\infty \frac{H_n}{n^3}&=&\frac{5}{4}\zeta (4)\\
\sum_{n=1}^\infty \frac{H_n}{n^4}&=&3 \zeta (5) - \zeta (2) \zeta (3) \, .
\end{eqnarray}

Since the functional equation of the Riemann zeta - function relates $\zeta (2k+1)$ to $\zeta '(-2k)$, this further convinces us of our conjecture. Such a relation exists, indeed! 

We state below the organization of the paper along with the main results and a sketch of the proof of the main theorem (Theorem 1). 

In Section \ref{sect:ramanuansum} we define the Ramanujan summation, state its properties and give a few examples that will be used in this paper. In Section \ref{sect:ramansumexpgenfn} we obtain the formula connecting $\zeta '(-k)$ to $\sum_{n\geq 1}^{\mathcal{R}}n^k H_n$. This section consists of the following results.
Let $\displaystyle{\psi(x) = \frac{\Gamma'(x)}{\Gamma (x)}}$ denote the digamma function. 

\begin{lemma} For $|z|<\pi $, we have
\begin{eqnarray}
\sum_{n\geq 1}^{\mathcal{R}}e^{nz}\log n&=&\frac{e^z}{e^{z}-1}%
\int_{0}^{1}\psi (x+1)e^{-zx}dx \nonumber \\
&&+\sum_{k=1}^{\infty} \frac{z^kH_k}{k!}\left(\frac{1-B_{k+1}}{k+1}\right) \nonumber \\
&&+\sum_{k= 1}^\infty \frac{z^{k}}{k!}%
\sum_{m=1}^{k+1}\frac{B_{m}}{m!}\frac{H_{k-m+1}}{(k-m+1)!} \, . \label{eq:enzlogn}
\end{eqnarray}
\end{lemma}

\begin{lemma}
For $|z|<\pi $
\begin{eqnarray}
\sum_{n\geq 1}^{\mathcal{R}}e^{-nz}H_{n} &=&\frac{1}{1-e^{-z}}\log  \left (\frac{z}{%
1-e^{-z}}\right )+\frac{1}{1-e^{-z}}\sum_{n= 1}^{\infty}\frac{%
(-1)^{n}}{n}\frac{z^{n}}{n!} \nonumber \\
&&+\gamma (\frac{1}{1-e^{-z}}-\frac{e^{-z}}{z})-\frac{1}{e^{z}-1}%
\int_{0}^{1}\psi (x+1)e^{-zx}dx \, . \label{eq:eznhn}
\end{eqnarray}
\end{lemma}
Note that the integral $\displaystyle{\int_{0}^{1}\psi (x+1)e^{-zx}dx}$ appears on the right hand side of both the equations (\ref{eq:enzlogn}) and (\ref{eq:eznhn}) and can be eliminated. This is the key idea of our main theorem which we state below.
\begin{theorem} For $|z|<\pi $
\be
\sum_{n\geq 1}^{\mathcal{R}}e^{nz}\log n+e^{z}\sum_{n\geq 1}^{\mathcal{R}%
}e^{-nz}H_{n}=\gamma \left(\frac{e^{z}}{1-e^{-z}}-\frac{1}{z}\right)\text{ ~mod~}(\mathbb{%
Q[[}z\mathbb{]])} \, . \label{eq:maintheorem}
\ee
\end{theorem}
Equating the coefficients of $z^k$,  we have 
\begin{corollary} If $k\geq 0$ then
\begin{equation}
\zeta ^{\prime }(-k)=\sum_{n\geq 1}^{\mathcal{R}}(1-n)^{k}H_{n}-\gamma \left(%
\frac{B_{k+1}(2)}{k+1}\right)\text{ mod~} \mathbb{Q} \, .
\end{equation}
When $k \ge 1$, this means that
\begin{equation}
\zeta ^{\prime }(-k)=\sum_{n\geq 1}^{\mathcal{R}}(1-n)^{k}H_{n}-\gamma \left(%
\frac{B_{k+1}}{k+1} +1\right)\text{ mod~} \mathbb{Q} \, .
\end{equation}
Equivalently for any integer $k\geq 1$ we have 
\begin{equation}
\sum_{n\geq 1}^{\mathcal{R}}n^{k}H_{n}=\gamma (\frac{1-B_{k+1}}{k+1})-\log (%
\sqrt{2\pi })+\sum_{m=1}^{k}\frac{k!}{m!(k-m)!}(-1)^{m}\zeta ^{\prime }(-m)%
\text{ mod ~}\mathbb{Q} \, , \label{eq:cor}
\end{equation}
and when $k=0$
\begin{equation}
\sum_{n\geq 1}^{\mathcal{R}}H_{n}=\frac{3}{2}\gamma -\log (\sqrt{2\pi })+%
\frac{1}{2} \, .
\end{equation}
\end{corollary}
The notations $\text{ ~mod~}(\mathbb{%
Q[[}z\mathbb{]])}$ and mod $\mathbb{Q}$ are used only for the simplicity of the statements of the results and these terms can be explicitly calculated as can be seen from the proof of the theorem. Further, we will express at the end of this section $\displaystyle{\sum_{n\geq 1}^{%
\mathcal{R}}n^{k}H_{n}}$ in terms of the analytic continuation of the normal
sum $\displaystyle{\sum_{n=1}^{\infty }\frac{H_{n}}{n^{s}}}$ near the
negative integers.

Next, in Section \ref{sect:interpretation}, we will show that the Laplace transforms and the Borel sums arise naturally when we consider the Ramanujan summation of exponential generating functions. Thus we interpret our results on the Ramanujan summation in classical terms as follows.
\begin{theorem} For $o < z < \pi $, we have
\begin{eqnarray}
\lefteqn{\sum_{k=0}^{\infty }\frac{z^{k}}{k!}\zeta ^{\prime }(-k)- \frac{\gamma}{z}%
+\sum_{n\geq 1}^{Borel}\frac{(-1)^{n+1}n!}{z^{n+1}}\zeta (n+1) =~~~~~~~~~~~~~} \nonumber \\
&&\frac{1}{%
e^{z}-1}\mathcal{L}\left (\frac{1}{x+1}\right )(z) 
-\sum_{k\ge 0} \frac{z^{k}}{(k+1)!}\frac{1}{k+1}- \nonumber \\
&&\sum_{k\geq 0}\frac{z^{k}}{k!}\sum_{m=1}^{k}\frac{B_{m}}{m!}\frac{%
H_{k-m+1}}{(k-m+1)!} 
+\sum_{k\geq 1}\frac{z^{k}}{k!}H_{k}\left (\frac{1-B_{k+1}}{k+1}\right ) \, ,
\end{eqnarray}
where $\mathcal{L}(f)$ denotes the Laplace transform of the function $f$. 
\end{theorem}
Also we show that an asymptotic formula obtained by Ramanujan for $\displaystyle{\sum_{n\geq 1} e^{-nz}\log (n)}$ is actually an equality. Section \ref{sect:higherharmonic} consists of a generalization of Lemma 2. Here we calculate $\displaystyle{\sum_{n\geq 1}^{\mathcal{R}}H_n^{(j)}~e^{-nz}}$ where $\displaystyle{H_n^{(j)}=\sum_{m=1}^n \frac{1}{m^j}}$.

\bigskip

\section{Ramanujan Summation} \label{sect:ramanuansum}  We refer to \cite{Candelpergher:Coppo:Delabaere} for the proofs of the results in this section. We give the definition of the Ramanujan summation in Section \ref{sect:def}, state its properties in Section \ref{sect:property} and furnish a few examples in Section \ref{sect:example}.

\subsection{Definition} \label{sect:def}
We will now define the Ramanujan summation of an infinite series. Let $f \in C^{\infty }(\mathbb{R})$.  We have, for $N\ge 1$, the Euler-MacLaurin formula 
\begin{equation}
f(0)=\int_{0}^{1}f(t)dt+\sum_{n=1}^{N}\partial ^{n-1}f(t)]_{0}^{1}\frac{B_{n}%
}{n!}+(-1)^{N+1}\int_{0}^{1}\partial ^{N}f(t)\frac{B_{N}(t)}{N!}dt \, .
\end{equation}
By a translation on $f$ we have for every integer $k$ 
\begin{equation}
f(k)=\int_{k}^{k+1}f(t)dt+\sum_{n=1}^{N}\partial ^{n-1}f(t)]_{k}^{k+1}\frac{%
B_{n}}{n!}+(-1)^{N+1}\int_{k}^{k+1}\partial ^{N}f(t)\frac{b_{N}(t)}{N!}dt \, , \label{eq:fk}
\end{equation}
where $b_{N}(t)=B_{N}(t-[t])$. 
So, if the integral $\displaystyle{\int_{1}^{\infty }\partial ^{N}f(t)\frac{b_{N}(t)}{N!}dt}$ is convergent, then we have
\begin{equation}
f(k)=R(k)-R(k+1) \, ,
\end{equation}
where 
\begin{equation}
R(x) =-\int_{1}^{x}f(t)dt-\sum_{k=1}^{N}\frac{B_{k}}{k!}%
\partial ^{k-1}f(x) +(-1)^{N+1}\int_{x}^{\infty}\partial ^{N}f(t)\frac{b_{N}(t)}{N!}dt \, .\label{eq:Rx} 
\end{equation}
Then a summation for $k=1,\cdots ,n-1$ gives 
\begin{equation}
\sum_{k=1}^{n-1}f(k)=R(1)-R(n) \, .
\end{equation}
The constant 
\begin{equation}
C_{f}=R(1)=-\sum_{k=1}^{N}\frac{B_{k}}{k!}\partial
^{k-1}f(1)+(-1)^{N+1}\int_{1}^{\infty }\partial ^{N}f(t)\frac{b_{N}(t)}{N!}dt 
\end{equation}
is the Ramanujan constant of the series $\sum_{n\geq 1}f(n)$. In Chapter
6 of his Notebooks Ramanujan says: ``The constant of a series has some
mysterious connection with the given infinite series and it is like the
centre of gravity of a body. Mysterious because we may substitute it for the
divergent series."

If all the derivatives of $f$ are sufficiently decreasing at infinity we can
verify, by integration by parts, that $R(1)$ does not depend of $N$ and that
the function $R$ satisfies the condition 
\begin{equation}
\lim_{x\rightarrow \infty }R(x)=-\int_{1}^{\infty }f(t)dt
\end{equation}
If we want to relax the hypothesis on the function $f$ we can try to define
the constant $C_{f}$ to be the value at $1$ of a solution $R$ of the
difference equation
\begin{equation}
R(x)-R(x+1)=f(x)\, .
\end{equation}
To obtain a unique value $R(1)$ we must find some conditions which gives us
a unique solution of this equation. We observe that for any integer $k$ 
\begin{equation}
\int_{k}^{k+1}f(x)dx=\int_{k}^{k+1}R(x)dx-\int_{k+1}^{k+2}R(x)dx \, .
\end{equation}
If $\lim_{x\rightarrow \infty }R(x)$ exists we have by a summation for $%
k\geq 1$ 
\begin{equation}
\int_{1}^{\infty }f(x)dx=\int_{1}^{2}R(x)dx-\lim_{x\rightarrow \infty }R(x) \, .
\end{equation}
Thus the condition $\lim_{x\rightarrow \infty }R(x)=-\int_{1}^{\infty
}f(x)dx$ is equivalent to the condition 
\begin{equation}
\int_{1}^{2}R(x)dx=0 \, .\label{eq:intR12}
\end{equation}

Unfortunately this condition does not give the uniqueness of the solution of the
difference equation as we can add to a solution any periodic function of period 
$1.$ To avoid this we ask that the function $R$ is analytic and of the order of
increase less than $2\pi .$

The following theorem provides us the existence of a unique solution for $R(x)$. See Section 3.1 of \cite{Candelpergher:Coppo:Delabaere}.
\begin{theorem}  Let $x\rightarrow a(x)$ be an analytic function on $P=\{x|\func{Re} (x)>0\}$ such
that there exist $\alpha <2\pi $ and $K>0$ with 
\begin{equation}
|a(x)|\leq Ke^{\alpha |x|}
\end{equation}
for all $x\in P.$ The equation 
\begin{equation}
R(x)-R(x+1)=a(x) \label{eq:difference}
\end{equation}
has only one solution $R_{a}$ analytic on $P=\{x|\func{Re} (x)>0\}$ so that
there exist $\beta <2\pi $ and $L>0$ with 
\begin{equation}
|R_{a}(x)|\leq Le^{\beta |x|}
\end{equation}
and with 
\begin{equation}
\int_{1}^{2}R_{a}(t)dt=0 \, .
\end{equation}
\end{theorem}
Now we are ready to give the definition of the Ramanujan summation of a series.
\begin{definition}
Let $x\rightarrow a(x)$ be an analytic function on $P=\{x|\func{Re} (x)>0\}$ such
that there exist $\alpha <\pi $ and $K>0$ with 
\begin{equation}
|a(x)|\leq Ke^{\alpha |x|}
\end{equation}
for all $x\in P.$ We call the Ramanujan summation of the series $\displaystyle{%
\sum_{n\geq 1}a(n)}$ the number 
\begin{equation}
\sum_{n\geq 1}^{\mathcal{R}}a(n)=R_{a}(1)\, .
\end{equation}
\end{definition}

We make the following remarks.

\noindent{\bf Remark 1.} The condition $\alpha <\pi $ is sufficient to obtain $a=b$ if $a(n)=b(n)$
for all $n\in \mathbb{N}$ and thus 
\begin{equation}
a(n)=b(n)\text{ for all }n\in \mathbb{N}\Rightarrow \sum_{n\geq 1}^{\mathcal{%
R}}a(n)=\sum_{n\geq 1}^{\mathcal{R}}b(n)\, .
\end{equation}

\noindent{\bf Remark 2.} We have 
\begin{equation}
R_{a}(1)-R_{a}(N)=\sum_{n<N}a(n)\, .
\end{equation}
If $\lim_{N\rightarrow \infty }R_{a}(N)$ exists, then the series $\displaystyle{\sum_{n=1}^\infty a(n)}$ is convergent and hence
\begin{equation}
\sum_{n\geq 1}^{\mathcal{R}}a(n)=\sum_{n=1}^{\infty
}a(n)+\lim_{N\rightarrow \infty }R_{a}(N)\, .
\end{equation}
Integrating the difference equation (\ref{eq:difference}) from $1$ to $N$ and using (\ref{eq:intR12}) we get
\begin{equation}
\int_{N}^{N+1}R(x)dx=-\int_{1}^{N}a(x)dx\, .
\end{equation}
So we obtain 
\begin{equation}
\sum_{n\geq 1}^{\mathcal{R}}a(n)=\sum_{n=1}^{\infty}a(n)-\int_{1}^{\infty
}a(x)dx \, .\label{eq:cgt}
\end{equation}

\subsection{Properties} \label{sect:property} 
We list below the properties of the Ramanujan summation. 

\begin{property} Linearity. 
\begin{equation}
\sum_{n\geq 1}^{\mathcal{R}}\lambda a(n)+\mu b(n)=\lambda \sum_{n\geq 1}^{%
\mathcal{R}}a(n)+\mu \sum_{n\geq 1}^{\mathcal{R}}b(n) \, .
\end{equation}
\end{property}

\begin{property} Translation. 
\begin{equation}
\sum_{n\geq 1}^{\mathcal{R}}a(n+N)=\sum_{n\geq 1}^{\mathcal{R}%
}a(n)-[a(1)+...+a(N)]+\int_{1}^{N+1}a(t)dt \, .
\end{equation}
\end{property} 
Note that this property is different from the axiom (C) Hardy states in Section 1.3, Chapter I of \cite{hardy}.
\begin{property} Derivation. 
\begin{equation}
R_{\partial ^{k}a}=\partial ^{k}R_{a}+\partial ^{k-1}a(1) \, , \label{eq:derivation}
\end{equation}
and
\begin{equation}
\sum_{n\geq 1}^{\mathcal{R}}\partial ^{k}a(n)=\partial ^{k}R_{a}(1)+\partial
^{k-1}a(1) \, .
\end{equation}
\end{property}

\begin{property} Taking $f=R_a$ and $k=1$ in (\ref{eq:fk}) and using (\ref{eq:derivation}) we get
\begin{equation}
\sum_{n\geq 1}^{\mathcal{R}}a(n)=-\sum_{n=1}^{N}\partial ^{n-1}a(1)\frac{%
B_{n}}{n!}+(-1)^{N+1}\int_{1}^{2}R_{\partial ^{N}a}(t)\frac{B_{N}(t-1)}{N!}dt \, .
\end{equation}
\end{property}

\begin{property} If $(z,x)\rightarrow a(z,x)$ is analytic on $D\times \{x|\func{Re} (x)>0\}$
and 
\begin{equation}
|a(z,x)|\leq Ce^{\alpha |x|}
\end{equation}
for some $\alpha <\pi $ and $C>0$, then $z\rightarrow \sum_{n\geq 1}^{%
\mathcal{R}}a(z,n)$ is analytic on $D$ and for all $z\in D$ 
\begin{equation}
\partial _{z}^{m}\sum_{n\geq 1}^{\mathcal{R}}a(z,n)=\sum_{n\geq 1}^{\mathcal{%
R}}\partial _{z}^{m}a(z,n) \, .
\end{equation}
In particular if $D$ is an open set containing $0$ and we have 
\begin{equation}
a(z,n)=\sum_{k\geq 0}a_{k}(n)z^{k} \, ,
\end{equation}
then for $z\in D$ 
\begin{equation}
\sum_{n\geq 1}^{\mathcal{R}}a(z,n)=\sum_{k\geq 0}\left(\sum_{n\geq 1}^{\mathcal{R}%
}a_{k}(n)\right)z^{k} \, . \label{eq:property5}
\end{equation}

\end{property}

\subsection{Examples} \label{sect:example} Now we give a few examples of the Ramanujan summation which would be used in the sequel.

\noindent{\bf Example 1.} If $a(x)=x^{-s}$ $(s\in \mathbb{C)}$, we have 
\begin{eqnarray}
R_{a}(x) &=&\zeta (x,s)-\frac{1}{s-1},~\text{ if }s\neq 1 \, ,\\
&=&-\frac{\Gamma ^{\prime }}{\Gamma }(x)=-\psi (x),~\text{ if }s=1 \, , \label{eq:Rpsi}
\end{eqnarray}
where
\be
\zeta (x,s) = \sum_{n=0}^\infty \frac{1}{(n+x)^s}\, .
\ee
Thus 
\begin{eqnarray}
\sum_{n\geq 1}^{\mathcal{R}}\frac{1}{n^{s}} &=&\zeta (s)-\frac{1}{s-1},~\text{
if }s\neq 1 \, ,\\
&=&\gamma,~ \text{ if }s=1 \, .
\end{eqnarray}
In particular 
\begin{equation}
\sum_{n\geq 1}^{\mathcal{R}}1=\frac{1}{2}
\end{equation}
and
\begin{equation}
\sum_{n\geq 1}^{%
\mathcal{R}}n^{k}=\frac{1-B_{k+1}}{k+1},~\text{ if }k=1,2,3,\cdots . \label{eq:sigmank}
\end{equation}
By derivation 
\begin{equation}
\sum_{n\geq 1}^{\mathcal{R}}\frac{\log n}{n^{s}}=-\zeta ^{\prime }(s)-\frac{1%
}{(s-1)^{2}},~\text{ if }s\neq 1 \label{eq:zeta'}
\end{equation}
and in particular 
\begin{equation}
\sum_{n\geq 1}^{\mathcal{R}}\log n=-\zeta ^{\prime }(0)-1=-1+\log  (\sqrt{%
2\pi }) \, . \label{eq:zeta'0}
\end{equation}

\noindent{\bf Example 2.} For $|z|<\pi $, let $a(x)=e^{xz}$. Then 
\begin{equation}
R_{a}(x)=\frac{e^{xz}}{1-e^{z}}+\frac{e^{z}}{z} 
\end{equation}
and 
\begin{equation}
\sum_{n\geq 1}^{\mathcal{R}}e^{nz}=\frac{e^{z}}{1-e^{z}}+\frac{e^{z}}{z} \, .\label{eq:ez}
\end{equation}
We observe also that $\displaystyle{z\rightarrow \sum_{n\geq 1}^{\mathcal{R}}e^{nz}}$ is
analytic for $|z|<\pi $, but the series $\displaystyle{\sum_{n\geq 1}e^{nz}}$ is convergent
for $\func{Re} (z)<0$ and so we have an analytic continuation of $\displaystyle{\sum_{n\geq 1}^{%
\mathcal{R}}e^{nz}}$ outside the disc $|z|<\pi $ if we define $\displaystyle{\sum_{n\geq
1}^{\mathcal{R}}e^{nz}}$ for $\func{Re} (z)<0$ by 
\begin{equation}
\sum_{n\geq 1}^{\mathcal{R}}e^{nz}=\sum_{n=1}^{\infty
}e^{nz}-\int_{1}^{\infty}e^{xz}dx=\frac{e^{z}}{1-e^{z}}+\frac{e^{z}}{z} \, .
\end{equation}
In the disc $|z|<\pi $ we have the relation 
\begin{equation}
\sum_{n\geq 1}^{\mathcal{R}}e^{-nz}=-e^{-z}\sum_{n\geq 1}^{\mathcal{R}%
}e^{nz}+\frac{1-e^{-z}}{z} \, . \label{eq:e-z}
\end{equation}

\section{Ramanujan Summation of the Exponential Generating Function} \label{sect:ramansumexpgenfn} In this section we will prove our main results. In Section \ref{sect:log} we will calculate $\displaystyle{\sum_{n\geq 1}^{\mathcal{R}} \frac{e^{-nz}}{n}}$. We will prove Lemma 1 in Section \ref{sect:enzlogn}, Lemma 2 in Section \ref{sect:harmonic} and Theorem 1 and its corollary in Section \ref{sect:proof}.

\subsection{Calculation of $\displaystyle{\sum_{n\geq 1}^{\mathcal{R}} \frac{e^{-nz}}{n}}$}\label{sect:log}
The Ramanujan sum ${\displaystyle \sum_{n\geq 1}^{\mathcal{R}}\frac{1}{n}e^{-nz}}$ is
defined in the disk $|z|<\pi $. For $0< \func{Re}(z) <\pi $ the series $\displaystyle{
\sum_{n=1}^\infty \frac{1}{n}e^{-nz}}$ is convergent. Thus using (\ref{eq:cgt}) we can write 
\begin{eqnarray}
\sum_{n\geq 1}^{\mathcal{R}}\frac{1}{n}e^{-nz} &=&\sum_{n=1}^{\infty}\frac{1}{n}e^{-nz}-\int_{1}^{\infty }\frac{1}{x}e^{-xz}dx \nonumber \\
&=&-\log  (1-e^{-z})+\func{Ei}(-z) \, ,
\end{eqnarray}
where the function $\func{Ei}$ is defined by \cite{gradshtein:ryzhik}
\begin{equation}
\func{Ei}(-z) =\left\{\begin{array}{ll} 
-\int_{z}^{\infty }\frac{1}{u}e^{-u}du,~& \text{ if }|\arg
(z)|<\pi \, , \\
&~\\
\frac{\func{Ei}(-z+i0)+\func{Ei}(-z-i0)}{2},~& \text{ if }z<0 \, . \end{array} \right .
\end{equation}
We have for $0<\func{Re}(z)<\pi $%
\begin{equation}
\func{Ei}(-z)=\gamma +\log z+\sum_{n=1}^{\infty}\frac{%
(-1)^{n}}{n}\frac{z^{n}}{n!} \, .
\end{equation}
Thus for $0<\func{Re}(z)<\pi $ 
\begin{equation}
\sum_{n\geq 1}^{\mathcal{R}}\frac{1}{n}e^{-nz}=\log  \left (\frac{z}{1-e^{-z}}%
\right )+\gamma +\sum_{n=1}^{\infty}\frac{(-1)^{n}}{n}\frac{z^{n}}{%
n!} \, .
\end{equation}
By analytic continuation this result remains true for $|z|<\pi .$ Thus we have proved
\begin{lemma} For $|z|<\pi $ we have 
\begin{equation}
\sum_{n\geq 1}^{\mathcal{R}}\frac{1}{n}e^{-nz}=\log \left (\frac{z}{1-e^{-z}}%
\right )+\gamma +\sum_{n=1}^{\infty}\frac{(-1)^{n}}{n}\frac{z^{n}}{%
n!} \, .
\end{equation}
\end{lemma}

\subsection{Proof of Lemma 1 - Calculation of $\displaystyle{\sum_{n\geq 1}^{\mathcal{R}} e^{nz}\log n}$}\label{sect:enzlogn}  Let $k$ be a non-negative integer and 
\begin{equation}
a_{k}(x)=\frac{x^{k}}{k!}(\log x -H_{k}) \, .
\end{equation}
It is easy to prove that \cite{roman}
\begin{equation}
\partial ^{m}a_k(x) =\left\{\begin{array}{ll}
\frac{x^{k-m}}{(k-m)!}(\log x -H_{k-m}),&\text{ if }0\leq
m\leq k  \, ,\\
~~\\
\frac{1}{x}, &\text{ if } m=k+1 \, .
\end{array} \right. 
\end{equation}
Thus using the linearity property, Property 4 of the Ramanujan summation and (\ref{eq:Rpsi}), we get
\begin{eqnarray}
\frac{1}{k!} \sum_{n\geq 1}^{\mathcal{R}}n^{k}\log n - \frac{H_k}{k!}\sum_{n\geq 1}^{%
\mathcal{R}}n^{k} &=&-\sum_{m=1}^{k+1}\frac{B_{m}}{m!}%
\partial ^{m-1}a_k(1)+(-1)^{k}\int_{1}^{2}R_{1/x}(t)\frac{B_{k+1}(t-1)}{(k+1)!}%
dt \nonumber \\
&=&\sum_{m=1}^{k}\frac{B_{m}}{m!}\frac{H_{k-m+1}}{(k-m+1)!}%
+(-1)^{k+1}\int_{1}^{2}\psi (t)\frac{B_{k+1}(t-1)}{(k+1)!}dt \, .\nonumber \\
&& ~~
\end{eqnarray}
The function 
\begin{equation}
f:z\rightarrow \sum_{n\geq 1}^{\mathcal{R}}e^{nz}\log n
\end{equation}
is analytic for $|z|<\pi $ and $\partial ^{k}f(0)=\sum_{n\geq 1}^{\mathcal{R}%
}n^{k}\log n $ and thus 
\begin{equation}
\sum_{n\geq 1}^{\mathcal{R}}e^{nz}\log n=\sum_{k=0}^{\infty}\frac{z^{k}}{%
k!}\sum_{n\geq 1}^{\mathcal{R}}n^{k}\log n \, .
\end{equation}
So we have
\begin{eqnarray}
\sum_{n\geq 1}^{\mathcal{R}}e^{nz}\log n &=& \sum_{k=0}^{\infty} \frac{z^kH_k}{k!}\sum_{n\geq 1}^{%
\mathcal{R}}n^{k}+\sum_{k\geq 0}\frac{z^{k}}{k!}%
\sum_{m=1}^{k}\frac{B_{m}}{m!}\frac{H_{k-m+1}}{(k-m+1)!}\nonumber \\
&&+\int_{1}^{2}\psi(t)\frac{1}{z}\sum_{k\geq 0}\frac{(-z)^{k+1}B_{k+1}(t-1)}{(k+1)!}dt \nonumber \\
&=&\sum_{k=1}^{\infty} \frac{z^kH_k}{k!}\left(\frac{1-B_{k+1}}{k+1}\right)+\sum_{k\geq 1}\frac{z^{k}}{k!}\sum_{m=1}^{k}\frac{B_{m}}{m!}\frac{%
H_{k-m+1}}{(k-m+1)!}\nonumber \\
&&+\frac{e^{z}}{e^z-1} \int_{0}^{1}\psi (t+1)e^{-zt} dt \, , \label{eq:enzln}
\end{eqnarray}
where we have used (\ref{eq:bernoulli}), (\ref{eq:sigmank}) and ${\displaystyle \int_1^2 \psi (t) dt = 0}$
This proves Lemma 1. 

\subsection{Proof of Lemma 2 - Calculation of $\displaystyle{\sum_{n\geq 1}^{\mathcal{R}} e^{-nz} H_n}$}\label{sect:harmonic}
By the linearity and the translation property of the Ramanujan summation, we can write for $%
|z|<\pi $ 
\begin{eqnarray}
\sum_{n\geq 1}^{\mathcal{R}}e^{-nz}H_{n} &=&e^{z}\sum_{n\geq 1}^{\mathcal{R}%
}e^{-(n+1)z}H_{n} \nonumber \\
&=&e^{z}\sum_{n\geq 1}^{\mathcal{R}}e^{-(n+1)z}H_{n+1}-e^{z}\sum_{n\geq 1}^{%
\mathcal{R}}e^{-(n+1)z}\frac{1}{n+1} \nonumber \\
&=&e^{z}\sum_{n\geq 1}^{\mathcal{R}}e^{-nz}H_{n}-1+e^{z}\int_{1}^{2}(%
\psi (x+1)+\gamma )e^{-zx}dx\nonumber \\
&& -e^{z}\sum_{n\geq 1}^{\mathcal{R}}\frac{%
e^{-(n+1)z}}{n+1} \, ,
\end{eqnarray}
where the last step follows from the fact that
\be
\psi(n+1) +\gamma =H_n \, . \label{eq:hn}
\ee
Now using the functional equation
\be
\psi(x+1) = \psi(x) +\frac{1}{x} \, , \label{eq:functionaleqnpsi}
\ee
we get
\begin{eqnarray}
\sum_{n\geq 1}^{\mathcal{R}}e^{-nz}H_{n} &=&\frac{-1}{1-e^{z}}+\frac{%
e^{z}}{1-e^{z}}\int_{1}^{2}\psi (x)e^{-zx}dx+\frac{e^{z}}{1-e^{z}}%
\int_{1}^{2}\frac{e^{-zx}}{x}dx-\gamma \frac{e^{-z}}{z} \nonumber \\
&&-\frac{e^{z}}{1-e^{z}}\sum_{n\geq 1}^{\mathcal{R}}\frac{e^{-(n+1)z}}{n+1} \, .
\end{eqnarray}
Using the translation property of the Ramanujan summation we write
\begin{equation}
\sum_{n\geq 1}^{\mathcal{R}}\frac{e^{-(n+1)z}}{n+1} = \sum_{n\geq 1}^{\mathcal{R}}\frac{e^{-nz}}{n}-e^{-z}+\int_1^2 \frac{e^{-tz}}{t}dt \, .
\end{equation}
Thus
\begin{equation}
\sum_{n\geq 1}^{\mathcal{R}}e^{-nz}H_{n}=\frac{e^z}{e^{z}-1}\sum_{n\geq 1}^{\mathcal{R}}\frac{1}{n}e^{-nz}-\gamma 
\frac{e^{-z}}{z}-\frac{e^z}{e^{z}-1}\int_{1}^{2}\psi (x)e^{-zx}dx \, . \label{eq:enzhn1}
\end{equation}
That is, for $|z|<\pi $
\begin{equation}
\sum_{n\geq 1}^{\mathcal{R}}e^{-nz}H_{n}=\frac{1}{1-e^{-z}}\sum_{n\geq 1}^{%
\mathcal{R}}\frac{1}{n}e^{-nz}-\gamma \frac{e^{-z}}{z}-\frac{1}{e^{z}-1}%
\int_{0}^{1}\psi (x+1)e^{-zx}dx \, .
\end{equation}
Lemma 2 now follows from Lemma 3.

\bigskip

\subsection{A Relation Between $\zeta '(-k)$ and $\sum_{n\geq 1}^{\mathcal{R}} n^k~H_{n}$ } \label{sect:proof} 
In this section we will prove Theorem 1 and Corollary 1. Note that the right hand side of both (\ref{eq:enzlogn}) and (\ref{eq:eznhn}) contain the term $\displaystyle{\int_0^1 \psi(x+1) e^{-xz} dx}$. Except for this term, all the other terms are known on the right hand side of both the equations. Surprisingly this unknown term can be eliminated and we get a simple relation between $\displaystyle{\sum_{n\geq 1}^{\mathcal{R}} e^{-nz} H_n}$ and $\displaystyle{\sum_{n\geq 1}^{\mathcal{R}} e^{nz} \log n}$. That is, we have for $|z|<\pi $,
\begin{eqnarray}
\sum_{n\geq 1}^{\mathcal{R}}e^{nz}\log n+e^{z}\sum_{n\geq 1}^{\mathcal{R}%
}e^{-nz}H_{n} &=&\gamma \left (\frac{e^{z}}{1-e^{-z}}-\frac{1}{z}\right ) 
+\sum_{k=1}^{\infty} \frac{z^kH_k}{k!}\left(\frac{1-B_{k+1}}{k+1}\right)\nonumber \\
&&+\sum_{k= 1}^{\infty} \frac{z^{k}}{k!}\sum_{m=1}^{k+1}\frac{B_{m}}{m!}\frac{%
H_{k-m+1}}{(k-m+1)!} \nonumber  \\
&&+\frac{e^{z}}{1-e^{-z}}\log \left (\frac{z}{1-e^{-z}}\right )\nonumber \\
&&+\frac{e^{z}}{1-e^{-z}}%
\sum_{n=1}^{\infty}\frac{(-1)^{n}}{n}\frac{z^{n}}{n!} \, . \label{eq:mainresult}
\end{eqnarray}
In other words, we have for $|z|<\pi $
\begin{equation}
\sum_{n\geq 1}^{\mathcal{R}}e^{nz}\log n+e^{z}\sum_{n\geq 1}^{\mathcal{R}%
}e^{-nz}H_{n}=\gamma (\frac{e^{z}}{1-e^{-z}}-\frac{1}{z})\text{ mod}(\mathbb{%
Q[[}z\mathbb{]])} \, .
\end{equation}
This proves Theorem 1.

\bigskip

Corollary 1 follows immediately by equating the coefficients of $z^k$. In particular, when $k=0$, we have from (\ref{eq:mainresult}) and (\ref{eq:zeta'})
\begin{equation}
\sum_{n\geq 1}^{\mathcal{R}}\log n+\sum_{n\geq 1}^{\mathcal{R}}H_{n}=\frac{3%
}{2}\gamma -\frac{1}{2} \, .
\end{equation}
Thus from (\ref{eq:zeta'0}) we have 
\begin{equation}
\sum_{n\geq 1}^{\mathcal{R}}H_{n}=\frac{3}{2}\gamma -\log (\sqrt{2\pi })+%
\frac{1}{2}
\end{equation}

\subsection{Interpretation of $\displaystyle{s\rightarrow \sum_{n\geq 1}^{\mathcal{R}}\frac{H_{n}}{n^{s}}}$} \label{sect:eulerharmonic} The function $\displaystyle{s\rightarrow \sum_{n\geq 1}^{\mathcal{R}}n^{-s}H_{n}}$ is an entire function and we have for $\func{Re}(s)>1$ the relation 
\begin{equation}
\sum_{n\geq 1}^{\mathcal{R}}n^{-s}H_{n}=h(s)-\int_{1}^{+\infty }x^{-s}(\psi
(x+1)+\gamma )dx \, ,
\end{equation}
where $h$ is the function defined for $\func{Re}(s)>1$ by  $\displaystyle{h(s)=\sum_{n\geq
1}^{\infty }n^{-s}H_{n}}$. This follows from (\ref{eq:cgt}) and (\ref{eq:hn}).

But we know the Taylor expansion of $x\rightarrow \psi (x+1)+\gamma $ and by
the interpolation formula of Ramanujan \cite{hardy1}, we obtain, if $1<Re(s)<2$, 
\begin{equation}
-\int_{0}^{\infty }x^{-s}(\psi (x+1)+\gamma )dx=\frac{\pi }{\sin (\pi s)}%
\zeta (s) \, .
\end{equation}
Thus, for $1<Re(s)<2$ 
\begin{equation}
h(s)=-\frac{\pi }{\sin (\pi s)}\zeta (s)-\int_{0}^{1}x^{-s}(\psi
(x+1)+\gamma )dx+\sum_{n\geq 1}^{\mathcal{R}}n^{-s}H_{n} \, .
\end{equation}
This shows that the function $h$  can be analytically extended for $Re(s)<1,$
with poles at the integers $0-1,-3,-5,...,1-2q,...$ with the expansions
\begin{equation}
h(s) = \left\{ \begin{array}{l}
\frac{1/2}{s}-\gamma -\zeta ^{\prime }(0)+\sum_{n\geq 1}^{\mathcal{R}%
}H_{n}+O(s)  \, ,\nonumber \\
~ \nonumber \\
\frac{\zeta (1-2q)}{s-(1-2q)}+\zeta ^{\prime
}(1-2q)-\int_{0}^{1}x^{2q-1}\psi (x)dx-\frac{\gamma }{2q} \nonumber \\
~ \nonumber \\
~~~+\sum_{n\geq 1}^{%
\mathcal{R}}n^{2q-1}H_{n}+O(s-(1-2q)) \, .
\end{array}
\right.
\end{equation}
For $s=-2q$, $q=1,2,...,$ we have $\zeta (-2q)=0$ and so the function $h$ is
analytic at $-2q$ and $\sum_{n\geq 1}^{\mathcal{R}}n^{2q}H_{n}$ is just
related to $h(-2q)$ by   
\begin{equation}
h(-2q)=-\zeta ^{\prime }(-2q)-\int_{0}^{1}x^{2q}\psi (x+1)dx-\frac{\gamma }{%
2q+1}+\sum_{n\geq 1}^{\mathcal{R}}n^{2q}H_{n} \, .
\end{equation}
The values of $h(-2q)$ are known in terms of Bernoulli numbers \cite{boyadzhiev:gadiyar:padma}.

\section{Ramanujan Summation of Exponential Generating Functions and the Borel Sums} \label{sect:interpretation} 
In this section we will express our results on the Ramanujan summation in classical terms. In Section \ref{sect:ramanujanlaplace} we will show that the Laplace transforms and the Borel sums appear naturally when we consider the Ramanujan summation of exponential generating functions. In Section \ref{sect:expborel} and Section \ref{sect:harmonicborel} we will interpret $\sum_{n\geq 1}^{\mathcal{R}} \frac{e^{-nz}}{n}$ and $\sum_{n\geq 1}^{\mathcal{R}}H_ne^{-nz}$ as Borel sums. In Section \ref{sect:ramanujanasymptotic} we prove that an asymptotic relation stated by Ramanujan for the sum $\displaystyle{\sum_{n=1}^\infty \log n ~e^{-nz}}$ is actually an equality. Finally we prove Theorem 2 in Section \ref{sect:zetaborel}.

\subsection{Ramanujan Summation, Laplace Transform and the Borel Sum} \label{sect:ramanujanlaplace}

Let $x\rightarrow f(x)$ be an analytic function on $P=\{x|\func{Re}(x)>0\}$ so
that for all $\varepsilon >0$ there exist $K>0$ with $|f(x)|\leq
Ke^{\varepsilon |x|}$ for all $x\in P.$

The Ramanujan summation of the series $\displaystyle{\sum_{n\geq 1}f(n)e^{-nz}}$ gives with
(\ref{eq:property5}), for $|z|<\pi $ 
\begin{equation}
\sum_{n\geq 1}^{\mathcal{R}}f(n)e^{-nz}=\sum_{k\geq 0}\frac{(-1)^{k}z^{k}}{k!%
}\sum_{n\geq 1}^{\mathcal{R}}n^{k}f(n) \, .
\end{equation}
Thus $\sum_{n\geq 1}^{\mathcal{R}}f(n)e^{-nz}$ is an exponential generating
function of the sums $\sum_{n\geq 1}^{\mathcal{R}}n^{k}f(n).$ These
generating functions are naturally related to Laplace transforms, because, if  
$0<z<\pi $, we have from (\ref{eq:cgt})  
\begin{equation}
\sum_{n\geq 1}^{\mathcal{R}}f(n)e^{-nz}=\sum_{n\geq 1}^{\infty
}f(n)e^{-nz}-e^{-z}\int_{0}^{\infty }e^{-xz}f(x+1)dx \, .
\end{equation}

\bigskip 

\begin{definition} Let $x\rightarrow f(x)$ be an analytic function on $P=\{x|\func{Re}(x)>0\}$ so
that for all $\varepsilon >0$ there exist $K>0$ with 
\begin{equation}
|f(x)|\leq Ke^{\varepsilon |x|}
\end{equation}
for all $x\in P.$ We define for $z>0$ the Laplace transform of $f$ by 
\begin{equation}
\mathcal{L(}f(x))(z)=\int_{0}^{\infty }e^{-xz}f(x)dx \, .
\end{equation}
\end{definition}
This transform is related to the Borel summation of divergent series of type 
$\sum_{n\geq 0}\frac{c_{n}}{z^{n+1}}.$

\begin{definition} 
A series $\sum_{n\geq 0}\frac{c_{n}}{z^{n+1}}$ is Borel summable for $z>0$
if the series $\sum_{n\geq 0}c_{n}\frac{x^{n}}{n!}$ has a radius of
convergence $R>0$ and if the function 
\begin{equation}
g(x)=\sum_{n=0}^{\infty }c_{n}\frac{x^{n}}{n!}
\end{equation}
has an analytic continuation along $\mathbb{R}_{+}$ with $\int_{0}^{\infty
}e^{-xz}g(x)dx$ convergent for $z>0.$ Then we define 
\begin{equation}
\sum_{n\geq 0}^{Borel}\frac{c_{n}}{z^{n+1}}=\int_{0}^{\infty }e^{-xz}g(x)dx \, .
\end{equation}
\end{definition}
\bigskip

\begin{lemma} Let $x\rightarrow f(x)$ be an analytic function on $P=\{x|\func{Re}(x)>0\}$ so
that for all $\varepsilon >0$ there exist $K>0$ with 
\begin{equation}
|f(x)|\leq Ke^{\varepsilon |x|}
\end{equation}
for all $x\in P$. Then
\begin{eqnarray}
\sum_{n\geq 1}^{\mathcal{R}}f(n)e^{-nz} &=&\sum_{n\geq 1}^{\infty
}f(n)e^{-nz}-e^{-z}\mathcal{L(}f(x+1))(z)  \\
&=&\sum_{n\geq 1}^{+\infty }f(n)e^{-nz}-e^{-z}\sum_{n\geq 0}^{Borel}\frac{%
\partial ^{n}f(1)}{z^{n+1}} \, .
\end{eqnarray}
\end{lemma}
\bigskip
\noindent{\bf Proof} From (\ref{eq:cgt}), we have for $0<z<\pi $, 
\begin{eqnarray}
\sum_{n\geq 1}^{\mathcal{R}}f(n)e^{-nz} &=&\sum_{n\geq 1}^{\infty
}f(n)e^{-nz}-\int_{1}^{\infty }e^{-xz}f(x)dx  \nonumber \\
&=&\sum_{n\geq 1}^{\infty }f(n)e^{-nz}-e^{-z}\mathcal{L(}f(x+1))(z) \, .
\end{eqnarray}
We have the Taylor expansion with the radius of convergence $\geq 1$ 
\begin{equation}
f(x+1)=\sum_{n\geq 0}\partial ^{n}f(1)\frac{x^{n}}{n!} \, .
\end{equation}
Then by definition 
\begin{equation}
\mathcal{L}(f(x+1))(z)=\sum_{n\geq 0}^{Borel}\frac{\partial ^{n}f(1)}{z^{n+1}} \, .
\end{equation}
Thus, for $0<z<\pi $, 
\begin{equation}
\sum_{n\geq 1}^{\mathcal{R}}f(n)e^{-nz}=\sum_{n\geq 1}^{+\infty
}f(n)e^{-nz}-e^{-z}\sum_{n\geq 0}^{Borel}\frac{\partial ^{n}f(1)}{z^{n+1}} \, .
\end{equation}
This proves Lemma 4.
\bigskip

\subsection{$\displaystyle{\sum_{n\geq 1}^{\mathcal{R}}\frac{e^{-nz}}{n}}$ as a Borel Sum} \label{sect:expborel} 
For $0<z<\pi $ we have from (\ref{eq:cgt})
\begin{eqnarray}
\sum_{n\geq 1}^{\mathcal{R}}\frac{1}{n}e^{-nz} &=&\sum_{n\geq 1}^{\mathcal{%
\infty }}\frac{1}{n}e^{-nz}-e^{-z}\mathcal{L}\left (\frac{1}{x+1}\right )(z) \nonumber \\
&=&-\log (1-e^{-z})-e^{-z}\mathcal{L}\left (\frac{1}{x+1}\right )(z) \, .
\end{eqnarray}
Thus, from Lemma 4,
\begin{equation}
\sum_{n\geq 1}^{\mathcal{R}}\frac{1}{n}e^{-nz}=-\log
(1-e^{-z})-e^{-z}\sum_{n\geq 0}^{Borel}\frac{(-1)^{n}n!}{z^{n+1}} \, .
\end{equation}
Note that this Borel sum is related to the $\func{Ei}$ function by the fact
that 
\begin{equation}
-e^{-z}\mathcal{L}\left (\frac{1}{x+1}\right )=\func{Ei}(-z) \, .
\end{equation}
\bigskip 

\subsection{$\displaystyle{\sum_{n\geq 1}^{\mathcal{R}}H_{n}e^{-nz}}$ as a Borel
Sum} \label{sect:harmonicborel}

For $0<z<\pi $ we have (again from (\ref{eq:cgt}))
\begin{eqnarray}
\sum_{n\geq 1}^{\mathcal{R}}H_{n}e^{-nz} &=&\sum_{n\geq 1}^{\mathcal{\infty 
}}H_{n}e^{-nz}-\int_{1}^{\infty }(\psi (x+1)+\gamma )e^{-xz}dx \nonumber \\
&=&\frac{-\log (1-e^{-z})}{1-e^{-z}}-\frac{\gamma e^{-z}}{z}%
-e^{-z}\int_{0}^{\infty }\psi (x+2)e^{-xz}dx \, .
\end{eqnarray}
The relation $\psi (x+2)=\psi (x+1)+\frac{1}{x+1}$ gives  
\begin{equation}
\sum_{n\geq 1}^{\mathcal{R}}H_{n}e^{-nz}=\frac{-\log (1-e^{-z})}{1-e^{-z}}-%
\frac{\gamma e^{-z}}{z}-e^{-z}\mathcal{L}(\psi (x+1))-e^{-z}\mathcal{L} \left (
\frac{1}{x+1}\right ) \, .
\end{equation}
Using the Taylor expansion 
\begin{equation}
\psi (x+1)=-\gamma +\sum_{n\geq 1}(-1)^{n+1}\zeta
(n+1)x^{n} \, ,
\end{equation}
this gives  
\begin{equation}
\sum_{n\geq 1}^{\mathcal{R}}H_{n}e^{-nz}=\frac{-\log (1-e^{-z})}{1-e^{-z}}%
-e^{-z}\sum_{n\geq 1}^{Borel}\frac{(-1)^{n+1}n!}{z^{n+1}}\zeta
(n+1)-e^{-z}\sum_{n\geq 0}^{Borel}\frac{(-1)^{n}n!}{z^{n+1}}
\end{equation}

\subsection{$\displaystyle{\sum_{n\geq 1}^{\infty } \log n~ e^{-nz}}$ and Generating
Function of $\protect\zeta ^{\prime }(-k)$} \label{sect:ramanujanasymptotic}

For $0<z<\pi $, we can write   
\begin{equation}
\sum_{n\geq 1}^{\mathcal{R}}\log n ~e^{-nz}=\sum_{n\geq 1}^{\infty }\log
n ~e^{-nz}-e^{-z}\mathcal{L}(\log (x+1))(z)\,.
\end{equation}
But $\mathcal{L}(\log (x+1))$ is simply related to the $\func{Ei}$ function.
We know from \cite{gradshtein:ryzhik} that 
\begin{eqnarray}
e^{-z}\mathcal{L}(\log (x+1))(z) &=&-\frac{1}{z}\func{Ei}(-z)  \notag \\
&=&\frac{-1}{z}(\gamma +\log z)+\sum_{k=0}^{\infty }\frac{(-1)^{k}z^{k}}{%
(k+1)!}\frac{1}{k+1}\,.
\end{eqnarray}
We know from (\ref{eq:zeta'}) that 
\begin{equation}
\sum_{n\geq 1}^{\mathcal{R}}e^{-nz}\log n=-\sum_{k=0}^{\infty }\frac{%
(-1)^{k}z^{k}}{k!}\zeta ^{\prime }(-k)-\sum_{k=0}^{\infty }\frac{%
(-1)^{k}z^{k}}{(k+1)!}\frac{1}{k+1}.
\end{equation}
$\,$
Hence we get 
\begin{equation}
\sum_{n\geq 1}^{\infty }e^{-nz}\log n=\frac{-1}{z}(\gamma +\log
z)-\sum_{k=0}^{\infty }\frac{(-1)^{k}z^{k}}{k!}\zeta ^{\prime }(-k)\,.
\end{equation}
Thus we see that Theorem 3.2 of Ramanujan in Chapter 15 of \cite{berndt1} is
not just an asymptotic expansion but an exact equality when $z\rightarrow 0+$
for the case $m=p=1$.

\subsection{$\displaystyle{\sum_{k=0}^{\infty }\frac{z^{k}}{k!}\protect\zeta ^{\prime }(-k)}
$ in Terms of Borel Sums}\label{sect:zetaborel}
In this section we will prove Theorem 2.
We have seen in (\ref{eq:enzln}) that for $|z|<\pi $ the sum $\sum_{n\geq 1}^{\mathcal{R}%
}e^{nz}\log n$ can be written in terms of the integral $\int_{0}^{1}\psi
(x+1)e^{-zx}dx.$ For $0<z<\pi $ there is a simple relation between $%
\int_{0}^{1}\psi (x+1)e^{-zx}dx$ and  the Laplace transform of $x\rightarrow
\psi (x+1)$. We have 
\begin{equation}
\mathcal{L(}\psi (x+1))(z)=e^{z}\int_{1}^{\infty }e^{-yz}\psi
(y)dy=e^{z}\int_{1}^{\infty }e^{-yz}\psi (y+1)dy-e^{z}\int_{1}^{\infty }%
\frac{e^{-yz}}{y}dy \, .
\end{equation}
Thus 
\begin{equation}
\mathcal{L}(\psi (x+1)=\frac{e^{z}}{e^{z}-1}\int_{0}^{1}e^{-yz}\psi (y+1)dy+%
\frac{1}{e^{z}-1}\mathcal{L}\left( \frac{1}{x+1}\right) (z)\,.
\end{equation}
For  recent developments on the Laplace transform of the digamma function, see for example, \cite{amdeberhan:moll} and \cite{glasser:manna}.
Now (\ref{eq:enzln}) becomes 
\begin{eqnarray}
\sum_{n\geq 1}^{\mathcal{R}}e^{nz}\log n &=&\mathcal{L(}\psi (x+1))(z)-\frac{%
1}{e^{z}-1}\mathcal{L}\left( \frac{1}{x+1}\right) (z)  \notag \\
&&+\sum_{k\geq 1}\frac{z^{k}}{k!}\sum_{m=1}^{k}\frac{B_{m}}{m!}\frac{%
H_{k-m+1}}{(k-m+1)!}+\sum_{k\geq 1}\frac{z^{k}}{k!}H_{k}\left( \frac{%
1-B_{k+1}}{k+1}\right) \,.  \notag \\
~&& ~
\end{eqnarray}
This can be translated in terms of Borel sums as follows.
\begin{eqnarray}
\sum_{n\geq 1}^{\mathcal{R}}e^{nz}\log n &=&-\frac{\gamma }{z}+\sum_{n\geq
1}^{Borel}\frac{(-1)^{n+1}n!}{z^{n+1}}\zeta (n+1)-\frac{1}{e^{z}-1}%
\sum_{n\geq 0}^{Borel}\frac{(-1)^{n}n!}{z^{n+1}}  \notag \\
&&+\sum_{k\geq 1}\frac{z^{k}}{k!}\sum_{m=1}^{k}\frac{B_{m}}{m!}\frac{%
H_{k-m+1}}{(k-m+1)!}+\sum_{k\geq 1}\frac{z^{k}}{k!}H_{k}\left( \frac{%
1-B_{k+1}}{k+1}\right) \,.   \notag \\
~&& ~ \label{eq:ramansumborelsum}
\end{eqnarray}

\bigskip 

But we know from (\ref{eq:zeta'}) that 
\begin{equation}
\sum_{n\geq 1}^{\mathcal{R}}e^{nz}\log (n)=-\sum_{k=0}^{+\infty }\frac{z^{k}%
}{k!}\zeta ^{\prime }(-k)-\sum_{k=0}^{\infty }\frac{z^{k}}{(k+1)!}\frac{1}{%
k+1}\,.  \label{eq:lognzeta'}
\end{equation}
Thus we obtain 
\begin{eqnarray}
\sum_{k=0}^{\infty }\frac{z^{k}}{k!}\zeta ^{\prime }(-k)-\frac{\gamma }{z}%
+\sum_{n\geq 1}^{Borel}\frac{(-1)^{n+1}n!}{z^{n+1}}\zeta (n+1) &=&\frac{1}{%
e^{z}-1}\mathcal{L}\left( \frac{1}{x+1}\right) (z)\notag \\
&& -\sum_{k=0}^{+\infty }%
\frac{z^{k}}{(k+1)!}\frac{1}{k+1}  \notag \\
&&-\sum_{k\geq 0}\frac{z^{k}}{k!}\sum_{m=1}^{k}\frac{B_{m}}{m!}\frac{%
H_{k-m+1}}{(k-m+1)!} \notag \\
&&+\sum_{k\geq 1}\frac{z^{k}}{k!}H_{k}\left( \frac{1-B_{k+1}}{k+1}\right) \,.
\end{eqnarray}
This proves Theorem 2.

\bigskip 

\noindent{\bf Remark 3. } 
We can write a very simple relation between these generating
functions

\begin{eqnarray}
\sum_{n\geq 1}^{\mathcal{R}}e^{nz}\log (n)+e^{z}\sum_{n\geq 1}^{\mathcal{R}%
}H_{n}e^{-nz} &=&\frac{e^{z}}{1-e^{-z}}\sum_{n\geq 1}^{\mathcal{R}}\frac{1}{n%
}e^{-nz}-\frac{\gamma }{z} \nonumber \\
&&+\sum_{k\geq 0}\frac{z^{k}}{k!}\sum_{m=1}^{k}\frac{B_{m}}{m!}\frac{%
H_{k-m+1}}{(k-m+1)!}\nonumber \\
&&+\sum_{k\geq 1}\frac{z^{k}}{k!}H_{k}(\frac{1-B_{k+1}}{k+1%
}) \, .
\end{eqnarray}
This is another form of Theorem 1.

\bigskip 

\section{Calculation of $\displaystyle{\sum_{n\geq 1}^{\mathcal{R}}H_n^{(j)}~e^{-nz}}$ } \label{sect:higherharmonic}
In this section we prove the following theorem. Let $Li_{j}(x) = \sum_{n=1}^\infty \frac{x^n}{n^j}$ for $0 <x <1$ denote the polylogarithm function.
\begin{theorem} For $0 < z< \pi$, we have
\begin{eqnarray}
\sum_{n\geq 1}^{\mathcal{R}}e^{-nz}H_{n}^{(j)} &=&\frac{1}{1-e^{-z}}%
Li_{j}(e^{-z}) -\frac{e^{-z}}{z}\zeta (j) \nonumber \\
&&+e^{-z}\sum_{m=1}^{j-1}z^{m-1}\frac{(j-m-1)!}{(j-1)!}\zeta \left(
j-m\right)\nonumber  \\
&&+\frac{e^{-z}}{1-e^{-z}}\sum_{m=1}^{j-1}z^{m-1}\frac{(j-m-1)!}{(j-1)!}%
(e^{-z}+(-1)^{m-1})\nonumber  \\
&&+\frac{(-1)^{j}z^{j-1}}{(j-1)!}e^{-z}\mathcal{L(}\psi (x+1))(z) 
+\frac{(-1)^{j}z^{j-1}}{(j-1)!}e^{-z}\mathcal{L(}\frac{1}{x+1})(z) \, . \nonumber \\
&&~ \label{eq:theorem4}
\end{eqnarray}
\end{theorem}

\noindent{\bf Proof} Using the functional equation (\ref{eq:functionaleqnpsi}) of $\psi $ we have
\begin{equation}
\partial ^{j-1}\psi (x+1)-\partial ^{j-1}\psi (x)=\partial ^{j-1}\frac{1}{x}=%
\frac{(-1)^{j-1}(j-1)!}{x^{j}} \, .
\end{equation}
That is,
\begin{equation}
\frac{(-1)^{j-1}}{(j-1)!}\partial ^{j-1}\psi (x+1)-\frac{(-1)^{j-1}}{(j-1)!}%
\partial ^{j-1}\psi (x)=\frac{1}{x^{j}} \, .
\end{equation}
Thus, 
\begin{equation}
H_{n}^{(j)}=\frac{(-1)^{j-1}}{(j-1)!}\partial ^{j-1}\psi (n+1)-\frac{%
(-1)^{j-1}}{(j-1)!}\partial ^{j-1}\psi (1) \, .
\end{equation}
But 
\begin{equation}
\partial ^{j-1}\psi (1)  =\left\{ \begin{array}{ll} (-1)^{j}(j-1)!\zeta (j),& \text{ ~if~ }j\geq 2  \\
-\gamma , \text{~ if ~}j=1 \, .\end{array} \right.
\end{equation}
which gives 
\begin{equation}
H_{n}^{(j)}=\frac{(-1)^{j-1}}{(j-1)!}\partial ^{j-1}\psi (n+1)+\zeta (j) \, .
\end{equation}
Here we have used the notation that $\zeta (1)=\gamma $.

For $0<z<\pi $ the series $\sum_{n\geq 1}H_{n}^{(j)}e^{-nz}$ is convergent and hence we can write  
\begin{eqnarray}
\sum_{n\geq 1}^{\mathcal{R}}H_{n}^{(j)}e^{-nz} &=&\sum_{n\geq 1}^{\mathcal{%
\infty }}H_{n}^{(j)}e^{-nz}-\int_{1}^{\infty }(\frac{(-1)^{j-1}}{(j-1)!}%
\partial ^{j-1}\psi (x+1)+\zeta (j))e^{-xz}dx \nonumber \\
&=&\sum_{n\geq 1}^{\mathcal{\infty }}H_{n}^{(j)}e^{-nz}-\frac{e^{-z}}{z}%
\zeta (j)-\frac{(-1)^{j-1}}{(j-1)!}\int_{1}^{\infty }\partial ^{j-1}\psi
(x+1)e^{-xz}dx \, . \nonumber \\
&&~~ \label{eq:hnjenz}
\end{eqnarray}
But we have 
\begin{equation}
\int_{0}^{1}e^{-xz}\partial ^{j-1}\psi (x+1)dx=\frac{e^{z}-1}{e^{z}}\mathcal{%
L(}\partial ^{j-1}\psi (x+1))-e^{-z}\mathcal{L(}\partial ^{j-1}\frac{1}{x+1}) \, ,
\end{equation}
and thus 
\begin{eqnarray}
\int_{1}^{\infty }\partial ^{j-1}\psi (x+1)e^{-xz}dx &=&\mathcal{L(}%
\partial ^{j-1}\psi (x+1))-\int_{0}^{1}\partial ^{j-1}\psi (x+1)e^{-xz}dx \nonumber \\
&=&(1-\frac{e^{z}-1}{e^{z}})\mathcal{L(}\partial ^{j-1}\psi (x+1))-e^{-z}%
\mathcal{L(}\partial ^{j-1}\frac{1}{x+1})\nonumber \\
&=&e^{-z}\mathcal{L(}\partial ^{j-1}\psi (x+1))-e^{-z}\mathcal{L(}\partial
^{j-1}\frac{1}{x+1}) \, .
\end{eqnarray}
Therefore, (\ref{eq:hnjenz}) becomes
\begin{eqnarray}
\sum_{n\geq 1}^{\mathcal{R}}H_{n}^{(j)}e^{-nz}  
&=&\sum_{n\geq 1}^{\mathcal{\infty }}H_{n}^{(j)}e^{-nz}-\frac{e^{-z}}{z}%
\zeta (j)-\frac{(-1)^{j-1}}{(j-1)!}e^{-z}\mathcal{L(}\partial ^{j-1}\psi (x+1)) \nonumber \\
&&+%
\frac{(-1)^{j-1}}{(j-1)!}e^{-z}\mathcal{L}\left (\partial ^{j-1}\frac{1}{x+1}\right ) \, .
\end{eqnarray}
That is,
\begin{eqnarray}
\sum_{n\geq 1}^{\mathcal{R}}H_{n}^{(j)}e^{-nz} &=&\frac{1}{1-e^{-z}}%
Li_{j}(e^{-z})-\frac{e^{-z}}{z}%
\zeta (j)-\frac{(-1)^{j-1}}{(j-1)!}e^{-z}\mathcal{L(}\partial ^{j-1}\psi (x+1)) \nonumber \\
&&+%
\frac{(-1)^{j-1}}{(j-1)!}e^{-z}\mathcal{L}\left (\partial ^{j-1}\frac{1}{x+1}\right ) \, . \label{eq:hnjenz1}
\end{eqnarray}
The Laplace transforms can be evaluated for $j\geq 2$ by integrating by parts repeatedly and thus we get
\begin{equation}
\mathcal{L(}\partial ^{j-1}f(x+1)) =-\sum_{n=2}^{j}z^{n-2}\partial ^{j-n}f(1)+z^{j-1}\mathcal{L(}f(x+1)) \, .
\end{equation}
Using this identity in (\ref{eq:hnjenz1}) for the function $f = \psi $, we get (\ref{eq:theorem4}). This proves Theorem 4.

\bigskip

\noindent \textbf{Remark 4. Relation Between $Li_{j}$ and the Bernoulli Numbers}

We have \cite{gessel}
\begin{equation}
\log \left (\frac{1-e^{-z}}{z}\right )=\sum_{n= 1}^\infty \frac{B_{n}}{n}\frac{z^{n}}{n!}
\end{equation}
and thus
\begin{equation}
Li_{1}(e^{-z})=-\log (1-e^{-z})=-\log z-\sum_{n= 1}^\infty \frac{B_{n}}{n}\frac{%
z^{n}}{n!} \, .
\end{equation}
Since $\partial _{z}Li_{2}(e^{-z})=-Li_{1}(e^{-z})$ we have 
\begin{equation}
Li_{2}(e^{-z})=C+z\log z-z+\sum_{n=11}^{\infty }\frac{B_{n}}{n}\frac{1}{n+1}%
\frac{z^{n+1}}{n!}\,,
\end{equation}
where $C=Li_{2}(1)=\zeta (2)\,.$ Again 
\begin{equation}
Li_{3}(e^{-z})=\zeta (3)-\zeta (2)z-\frac{1}{2}z^{2}\log z+\frac{3}{4}%
z^{2}-\sum_{n=11}^{\infty }\frac{B_{n}}{n}\frac{1}{n+1}\frac{1}{n+2}\frac{%
z^{n+2}}{n!}
\end{equation}
and in general for $j\geq 1$ we have 
\begin{eqnarray}
Li_{j}(e^{-z}) &=&\zeta (j)-\zeta (j-1)z+...+\zeta (2)\frac{(-1)^{j-1}z^{j-2}%
}{(j-2)!}  \notag \\
&&+\frac{(-1)^{j-1}z^{j-1}}{(j-1)!}(\log z-H_{j-1})  \notag \\
&&+(-1)^{j-1}\sum_{n=1}^{\infty }\frac{B_{n}}{n}\frac{1}{n+1}\frac{1}{n+2}...%
\frac{1}{n+j-1}\frac{z^{n+j-1}}{n!}\,.
\end{eqnarray}
(Note that $\partial ^{k}(\frac{z^{k}}{k!}(\log z-H_{k}))=\log z$.)

Thus we obtain 
\begin{eqnarray}
\sum_{n\geq 1}^{\mathcal{R}}H_{n}^{(j)} &=&\lim_{z\rightarrow 0}\frac{1}{%
1-e^{-z}}Li_{j}(e^{-z})-\frac{e^{-z}}{z}\zeta (j)  \notag \\
&&+\frac{(-1)^{j-1}}{(j-1)!}\partial ^{j-2}\psi (1)-\frac{(-1)^{j-1}}{(j-1)!}%
\partial ^{j-2}\left( \frac{1}{x}\right) (1)  \notag \\
&=&\lim_{z\rightarrow 0}\left( \frac{1}{1-e^{-z}}-\frac{e^{-z}}{z}\right)
\zeta (j)-\frac{z}{1-e^{-z}}\zeta (j-1) \nonumber \\
&=&\frac{3}{2}\zeta (j)-\zeta (j-1)\frac{j-2}{j-1}+1\,.
\end{eqnarray}

\section{Conclusion} This paper has two ingredients, the first being the calculus of the Ramanujan summation developed by Candelpergher, Coppo and Delabaere and the second the use of generating functions for divergent series. Using these two tools, we have obtained a formula connecting the derivatives of the Riemann zeta - function at the negative integers to the Ramanujan sum of divergent Euler sums. We have shown that the Ramanujan sum $\sum_{n\geq 1}^{\mathcal{R}%
} \log n~e^{-nz}$ gives naturally the exponential generating function of $%
\zeta ^{\prime }(-k)$ where $k=0,1,2,...$, and surprisingly this sum is
related to the Laplace transform $\mathcal{L(}\psi (x+1))$. We have also obtained an
explicit expression for the strange function $\displaystyle{\sum_{k=0}^{+\infty }\frac{z^{k}}{k!}\zeta ^{\prime }(-k)}-\frac{\gamma }{z}%
+\displaystyle{\sum_{n\geq 1}^{Borel}\frac{(-1)^{n+1}n!}{z^{n+1}}\zeta (n+1)}$. Further we have proved that an asymptotic formula for $\sum_{n=1}^{\infty
} \log n~e^{-nz}$ obtained by Ramanujan is actually an equality. Since the harmonic numbers $H_{n}$ are naturally interpolated by the function 
$\psi (x+1)+\gamma ,$ the Ramanujan sum $\displaystyle{\sum_{n\geq 1}^{\mathcal{R}%
}H_{n}e^{-nz}}$ is also related to the Laplace transform $\displaystyle{\mathcal{L(}\psi
(x+1))}$ and thus this gives an explicit expression for the values $\displaystyle{%
\sum_{n\geq 1}^{\mathcal{R}}n^{k}H_{n}}$ which are related to the analytic
continuation of the function $\sum_{n\geq 1}\frac{H_{n}}{n^{s}}$ near $s = -k$.


\begin{thebibliography}{15}
%
\bibitem{amdeberhan:moll} Amdeberhan, T, Espinosa, O., Moll, V. H.: The Laplace transform of the digamma function: an integral due to Glasser, Manna and Oloa. Proc. Amer. Math. Soc. {\bf 136} , 3211-3221 (2008).
%
\bibitem{apostol} Apostol, T. M.: Introduction to Analytic Number Theory. Springer International Student Edition (1980)
%
\bibitem{berndt} Berndt, B: Ramanujan's Notebooks, Part I. Springer (1985)
%
\bibitem{berndt1} Berndt, B: Ramanujan's Notebooks, Part II. Springer (1985)
%
\bibitem{boyadzhiev:gadiyar:padma} Boyadzhiev, K. N., Gadiyar, H. G., Padma, R.: A note on the values of an Euler sum at negative integers and relation to a convolution of Bernoulli numbers. Bull. Korean Math. Soc. {\bf 45}(2) 277-283 (2008)
%
\bibitem{Candelpergher:Coppo:Delabaere} Candelpergher, B., Coppo, M. A., Delabaere, E.: La sommation de Ramanujan. L'Enseignement Math\'{e}matique, {\bf 43}, 93-132 (1997) 
%
\bibitem{Flajolet:Salvy} Flajolet, P., Salvy, B.: Euler sums and contour integral representations. Experimental Mathematics {\bf 7}(1) 15-35 (1998)
%
\bibitem{gadiyar:padma} Gadiyar, H. G., Padma, R.: A comment on Matiyasevich's identity \#0102 with Bernoulli numbers, \url{http://arxiv.org/abs/math/0608675v2}
%
\bibitem{gessel} I. M. Gessel, \emph{On Miki's identities for Bernoulli numbers}, J. Number Theory, {\bf 110}  75-82 (2005).
%
\bibitem{glasser:manna} Glasser M, L, Manna, D.: On the Laplace transform of the Psi function. ``Tapas in Experimental Mathematics'' (T. Amdeberhan and V. Moll, eds.), Contemporary Mathematics {\bf 457} 193 - 202 (2008) 
%
\bibitem{gradshtein:ryzhik} Gradshteyn, I.S., Ryzhik, I.M.: Table of Integrals, Series, and Products.
Alan Jeffrey (ed.), Fifth edition, Acadamic Press (1994)
%
\bibitem{hardy1} Hardy, G. H.: Ramanujan: Twelve lectures on subjects suggested by his life and work. Chelsea Publishing Company (1940)
%
\bibitem{hardy} Hardy, G. H.: Divergent Series. Clarendon Press, Oxford (1949)
%
\bibitem{roman} Roman, S.: The logarithmic binomial formula. Amer. Math. Monthly {\bf 99}(7) 641-648 (1992)
%
\end{thebibliography}
\end{document}